# Cochain algebras of mapping spaces and finite group actions.

Frédéric Patras[(1)], Jean-Claude Thomas[(2)]

**Abstract.** The purpose of the present article is threefold. First of all, we rebuild the whole theory of cosimplicial models of mapping spaces by using systematically Kan adjunction techniques. Secondly, given two topological spaces $X$ and $Y$, we construct a cochain algebra which is quasi-isomorphic (as an algebra) to the singular cochain algebra of the mapping space $Y^X$. Here, $X$ has to be homotopy equivalent to the geometric realization a finite simplicial set and of dimension less or equal to the connectivity of $Y$. At last, we apply these results to the study of finite group actions on mapping spaces that are induced by an action on the source.

(1) CNRS UMR 6621 - Université de Nice, Mathématiques, Parc Valrose, 06108 Nice cedex 2, France. *patras@math.unice.fr* Tel. 00 33 (0) 492076262, Fax. 00 33 (0) 493517974.

(2) Université d'Angers & CNRS UMR 6093. Département de mathématiques 2, Bd Lavoisier 49045 Angers France. *jean-claude.thomas@univ-angers.fr*

# Introduction.

The methods of study of the cohomology of mapping spaces originate in two seminal papers of D.W. Anderson [1, 2]. The subject, besides its own interest, is related to various classical topics such as the cohomology of vector fields on a manifold (Gelfand-Fuchs), or the cohomology of configuration spaces. These connexions are explained by Bott-Segal and Bendersky-Gitler in [6, 5]. The same authors have developed tools for computing the cohomology of mapping spaces. They have established results that had been announced in Anderson's papers but had remained without a published proof. In particular, Bott and Segal have constructed a cochain complex computing the cohomology of a mapping space in the particular case when the source of the mapping space is provided with an open contractible covering (some finiteness and connectivity conditions are required, see section 2 of the present article for details). The simplices of the simplicial set associated to the covering (the nerve of the corresponding poset) are the building blocks of this cochain complex.

A cochain complex can be associated by the same process as in [6] to any simplicial set and any topological space. One would then expect the results of Bott-Segal to hold more generally for the corresponding mapping space (from the topological realization of the simplicial set to the topological space). That is, this cochain complex should compute the cohomology of the mapping space. A proof of this more general result along the lines of [6] has been announced in [5], but the proof by induction that is given in this article is not satisfactory : for example it does not to apply when the simplicial set is the usual simplicial model of the $n$-sphere.

Motivated by these observations, the purpose of the present article is threefold. First of all, we show how to rebuild the whole theory of cosimplicial models of mapping spaces by using systematically Kan adjunction techniques. This approach enlightens the theory but also provides new computational tools. Secondly, following the ideas of Anderson et al., given a simplicial set $K$ and a topological space $Y$, we construct a cochain algebra. We show that this cochain algebra is quasi-isomorphic (as an algebra) to the cochain algebra of singular cochains on the mapping space $Y^{|K|}$ (Thm. 2). Here, $K$ has to be a simplicial finite set weakly homotopy equivalent to a finite simplicial set and of dimension less or equal to the connectivity of $Y$. At last, we apply these results to the study of group actions on mapping spaces that are induced by a simplicial action on the source, reducing this study to representation theoretic arguments in a suitable category of functors. In particular, we show how the decomposition into isotypic components of the cohomology of the mapping space follows from the corresponding decomposition of the homology of the source, together with the multiplicative structure of the Grothendieck ring of the group.

Our results also hold *mutatis mutandis* when simplicial sets, spaces and mapping spaces are replaced by pointed simplicial sets, pointed spaces and pointed mapping spaces. In that case, the dual problem, that is the modelization and computation of the homology of pointed mapping spaces, has been studied by Bousfield by different means [7, 4.2]. The results of Bousfield follow from general theorems on the realization of cosimplicial spaces and on the (strong) convergence of the associated spectral sequences. When dealing with (pointed or not) mapping spaces, the (Kan) adjunction techniques we develop seem better suited to algebraic computations (modelization of products, group actions...) than these techniques.

**Content.**

# 1 Cosimplicial mapping spaces.

In this section we prove that cosimplicial mapping spaces can be thought of as right adjoints to a suitable tensor product functor -a very natural fact from the categorical point of view on mapping spaces. As a corollary of the unpointed case we recover a result of Bott-Segal, [6, Prop. 5.1]. In this section the category of topological spaces can be replaced by the category of simplicial sets.

## 1.1 Classical notations.

If **C** and **D** are two categories, we write:

- $A \in \mathbf{C}$ when $A$ is an object of **C**,

- $f \in \mathbf{C}(A, B)$ when $f : A \to B$ is a morphism of **C**,

- $\mathbf{C^{op}}$ for the category opposite to **C**,

- $F \in \mathbf{D^C}$ (resp. $F \in \mathbf{D^{C^{op}}}$) when $F$ is a covariant functor (resp. a contravariant functor) from **C** to **D**),

- **Set** for the category of sets,

- **Top** for the category of topological spaces (compactly generated). From now on, we call space a topological space and we write, as usual, $Y^X$ for $\mathbf{Top}(X, Y)$.

- $\mathbf{C}_*$ for the "pointed" category corresponding to **C**, whenever it makes sense. For example, $\mathbf{Set}_*$, (resp. $\mathbf{Top}_*$) is the category of pointed spaces (resp. of pointed spaces),

- $\mathbf{\Delta}$ is the category of finite ordered sets with objects $[n] := \{0, 1, ..., n\}$ and non decreasing maps (i.e. $i < j \Rightarrow f(i) \leq f(j)$).

- $\mathcal{C}os\mathbf{C} = \mathbf{C^{\Delta}}$ (resp. $\mathcal{S}\mathbf{C} = \mathbf{C^{\Delta^{op}}}$) the category of cosimplicial objects in **C** (resp. the category of simplicial objects in **C**).

- A cosimplicial object $\underline{L} : \mathbf{\Delta} \to \mathbf{C}$ (resp. a simplicial object $K : \mathbf{\Delta^{op}} \to \mathbf{C}$) is also written $(\underline{L}[n], d^i, s^i)$ (resp. $(K_n, d_n, s_n)$) where $\underline{L}[n] = \underline{L}([n])$ (resp. $K_n = K([n])$ ) and where $d^i, s^i$ (resp. $d_i, s_i$) stand for the (co)face and (co)degeneracy morphisms. A simplicial finite set is a simplicial set $K$ such that each component $K_n$ is finite. We say that a simplicial set is *finite* if it has a finite number of non-degenerate simplices.

## 1.2 Geometric realization of a cosimplicial space.

To each simplicial set $K$ is associated a space, written $|K|$ and called the topological realization of $K$ [15]. In the same way, each cosimplicial space has a geometric realization that is defined as follows. Consider the remarkable cosimplicial space $\underline{\Delta} \in \mathcal{C}os\mathbf{Top}$ such that $\underline{\Delta}[p] := \Delta^p$, the standard $p$-simplex in $\mathbf{R}^{p+1}$ ($\Delta^p = \{(t_0, ..., t_p) | \ 0 \leq t_i \leq 1, \ \sum_{i=0}^{p} t_i = 1\}$). The cosimplicial structure on $\underline{\Delta}$ is given as follows: if $\phi \in \mathbf{\Delta}([m], [n])$ then $\underline{\Delta}(\phi) \in \mathbf{Top}(\Delta^m, \Delta^n)$, $\underline{\Delta}(\phi)(t_0, ..., t_m) := (t'_0, ..., t'_n)$, where $t'_i := \sum_{\phi(j)=i} t_j$.
The geometric realization of a cosimplicial space [8], $\underline{X} \in \mathcal{C}os\mathbf{Top}$, is the space:
$$||\underline{X}|| := \mathcal{C}os\mathbf{Top}(\underline{\Delta}, \underline{X}).$$

The topology on $||\underline{X}||$ is induced by the canonical inclusion: $\mathcal{C}os\mathbf{Top}(\underline{\Delta}, \underline{X}) \subset \prod_{n \geq 0} \mathbf{Top}(\Delta^n, \underline{X}[n])$.

## 1.3 Left Kan extension and mixed tensor product under $\mathbf{\Delta}$.

Recall that, for a given category $\mathbf{C}$, the Yoneda embedding [4, 14], is the functor:

$$\mathcal{Y}_\mathbf{C} : \mathbf{C} \hookrightarrow \mathbf{Set}^{\mathbf{C}^{\mathrm{op}}} \begin{cases} C \mapsto \mathbf{C}(-, C), \\ (C \xrightarrow{f} C') \mapsto (\mathbf{C}(-, C) \xrightarrow{f_*} \mathbf{C}(-, C'), \varphi \mapsto f_*(\varphi) = f \circ \varphi). \end{cases}$$

We write $\mathcal{Y}$ for the Yoneda embedding: $\mathbf{\Delta} \hookrightarrow \mathbf{Set}^{\mathbf{\Delta}^{\mathrm{op}}} = \mathcal{S}\mathbf{Set}$. Thus, $\mathcal{Y}$ is a cosimplicial-simplicial set.

Since the category of topological spaces is cocomplete, and since $\mathbf{\Delta}$ is small, for any $\underline{Z} \in \mathcal{C}\mathcal{O}\mathcal{S}\mathbf{Top}$, there exists a left Kan extension of $\underline{Z}$ along $\mathcal{Y}$, [4, Sect. 3.7]. This is a functor : $\mathcal{S}\mathbf{Set} \to \mathbf{Top}$, usually denoted $Lan_\mathcal{Y} \underline{Z}$ but that we prefer to write using the simpler algebraic notation $\mathcal{Y} \otimes \underline{Z}$. Moreover the following diagram commutes up to homeomorphism since $\mathcal{Y}$ is full and faithful [4, Prop. 3.7.3].

$$
\begin{array}{ccc}
\mathbf{\Delta} & \xrightarrow{\mathcal{Y}} & \mathcal{S}\mathbf{Set} \\
{\scriptstyle \underline{Z}} \downarrow & \swarrow {\scriptstyle \mathcal{Y} \otimes \underline{Z}} & \\
\mathbf{Top} & &
\end{array}
$$

By analogy with the case of the tensor product of set-valued functors, for any $K \in \mathcal{S}\mathbf{Set}$ we put

$$(\mathcal{Y} \otimes \underline{Z})(K) := K \otimes_{\mathbf{\Delta}} \underline{Z}$$

and call this space *the tensor product under* $\mathbf{\Delta}$ of the simplicial set $K$ and of the cosimplicial space $\underline{Z}$.

This space $K \otimes_{\mathbf{\Delta}} \underline{Z}$ can also be described as a "classical" tensor product (defined by generators and relations). Let $\mathbf{C}$ and $\mathbf{D}$ be the categories defined as follows. The objects of $\mathbf{C}$ are the pairs $([n], \alpha)$, where $\alpha$ is a natural transformation from $\mathbf{\Delta}(-, [n])$ to $K$. An arrow of $\mathbf{C}$, $f : ([n], \alpha) \longrightarrow ([m], \beta)$ is a morphism $f \in \mathbf{\Delta}([n], [m])$ such that $\beta \circ \mathbf{\Delta}(-, f) = \alpha$.

The object of $\mathbf{D}$ are the simplices of $K$ and a morphism $\overline{f} : \sigma \to \tau$ with $\sigma \in K_n$ and $\tau \in K_m$ is a morphism $f \in \mathbf{\Delta}([n], [m])$ such that $K(f)(\tau) = \sigma$.

By the Yoneda lemma, $\mathbf{C}$ identifies with the category $\mathbf{D}$, up to a canonical isomorphism. Then, we have (see [4, Sect. 3.8]):

$$K \otimes_{\mathbf{\Delta}} \underline{Z} \cong colim_{([n],\alpha) \in \mathbf{C}} \; \underline{Z} \circ \mathcal{U}([n], \alpha) \cong \frac{\coprod_n \coprod_{\sigma \in K_n} \underline{Z}[n]}{\approx} \cong \frac{\coprod_n (K_n \times \underline{Z}[n])}{\approx},$$

where $\mathcal{U}$ denotes the forgetful functor $\mathbf{C} \longrightarrow \mathbf{\Delta}$, $\mathcal{U}([n], \alpha) := [n]$ and the equivalence relation $\approx$ is defined by:

$(\sigma, z) \approx (\tau, y)$ if there exists $f \in \mathbf{D}(\sigma, \tau)$ such that $y = \underline{Z}(f)(z)$.

For example, if $\underline{Z} = \underline{\Delta}$, $K \otimes_{\mathbf{\Delta}} \underline{\Delta}$ is the usual topological realization of $K$, see e.g. [15].

## 1.4 The pointed case.

The same constructions hold with pointed categories, up to some details. There is still a "pointed Yoneda embedding":

$$\mathcal{Y}_* : \mathbf{\Delta} \hookrightarrow \mathbf{Sets}_*^{\mathbf{\Delta}^{op}},$$

which is defined by requiring that the constant map $\tilde{0} \in \mathbf{\Delta}([m],[n])$, $\tilde{0}(i) = 0$, $i \in [n]$, is the "point" in $\mathbf{\Delta}([m],[n])$.

Using the same process as in the non-pointed case, we define, for any cosimplicial space $\underline{Z}$, $\mathcal{Y}_* \otimes \underline{Z} \in \mathbf{Top}^{\mathbf{Sets}_*^{\mathbf{\Delta}^{op}}}$ and for any simplicial pointed set $K_*$ the corresponding space $K_* \otimes_{\mathbf{\Delta}} \underline{Z}$.

Moreover, if $K_*$ is reduced (i.e. has a unique 0-simplex $*$), a direct inspection of the colimit formulas defining $K_* \otimes_{\mathbf{\Delta}} \underline{Z}$ shows that $K_* \otimes_{\mathbf{\Delta}} \underline{Z} \cong K \otimes_{\mathbf{\Delta}} \underline{Z}$, where $K$ is $K_*$ viewed as a simplicial set. This follows from the fact that, if $K_*$ is reduced, any natural transformation from $\mathbf{\Delta}(-,[n])$ to $K$ is automatically a map of pointed simplicial sets. In particular, if $K_*$ is reduced, $K_* \otimes_{\mathbf{\Delta}} \underline{\Delta}$ is isomorphic to the topological realization of $K_*$. Notice also that, if $K_*$ is reduced and $\underline{Z}$ is also reduced (in the sense that $\underline{Z}[0]$ is a point), then $K_* \otimes_{\mathbf{\Delta}} \underline{Z}$ is naturally pointed by the class of the unique element in $K_0 \times \underline{Z}[0]$.

## 1.5 The adjunction theorem for cosimplicial mapping spaces.

We introduce the *Hochschild-Yoneda functor* (resp. *the pointed Hochschild-Yoneda functor*) as the functor:

$$\mathcal{H} : \mathbf{Top} \to \mathbf{Top}^{\mathbf{Set}^{op}}, \quad (\text{resp.} \mathcal{H}_* : \mathbf{Top}_* \to \mathbf{Top}_*^{\mathbf{Set}_*^{op}})$$

defined by: $\mathcal{H}(X) := \mathbf{Set}(-, X)$ (resp. $\mathcal{H}_*(X_*) := \mathbf{Set}_*(-, X_*)$). The topology on $\mathbf{Set}(E, X) \cong X^E$ and $\mathbf{Set}_*(E_*, X_*) \cong X_*^{E_* - \{*\}}$ is the usual one (the product topology).

The terminology "Hochschild-Yoneda functor" is motivated by the observation that the definition of $\mathcal{H}$ mixes a kind of Yoneda-embedding process together with a topological generalization of the algebraic construction leading to the definition of Hochschild homology [21].

Let $K$ be a simplicial set and $X$ a space (resp. $K_*$ a simplicial pointed set and $X_*$ a pointed space). We call *cosimplicial mapping space from $K$ to $X$ (resp. from $K_*$ to $X_*$)* and write $X^K$ (resp. $X_*^{K_*}$) for the cosimplicial space:

$$X^K := \mathcal{H}(X) \circ K, \quad (resp. \ X_*^{K_*} : \mathcal{H}_*(X_*) \circ K_*).$$

In particular, if $K$ is finite $X^K[n] = \mathbf{Set}(K_n, X) = X^{\times \#K_n}$, where $\#K_n$ stands for the number of elements in $K_n$ and $X^{\times l}$ for the cartesian product of $l$ copies of $X$.

**Theorem 1** *(Adjunction theorem for cosimplicial mapping spaces.) Let $K$ be a simplicial set, $X$ a space and $\underline{Z}$ a cosimplicial space. There exists a natural homeomorphism*

$$\mathcal{C}os\mathbf{Top}(\underline{Z}, X^K) \cong \mathbf{Top}(K \otimes_{\boldsymbol{\Delta}} \underline{Z}, X).$$

*Let $K_*$ be a reduced simplicial set, $X_*$ a pointed space and $\underline{Z}$ a reduced cosimplicial space. There exists a natural homeomorphism*

$$\mathcal{C}os\mathbf{Top}(\underline{Z}, X_*^{K_*}) \cong \mathbf{Top}_*(K_* \otimes_{\boldsymbol{\Delta}} \underline{Z}, X_*).$$

These isomorphisms can be viewed as isomorphisms in **Top** by viewing $\mathcal{C}os\mathbf{Top}(\underline{Z}, X^K)$ (resp. $\mathcal{COS}\mathbf{Top}(\underline{Z}, X_*^{K_*})$) as a topological space for the embedding:

$$\mathcal{C}os\mathbf{Top}(\underline{Z}, X^K) \subset \prod_n \mathbf{Top}(\underline{Z}[n], X^K[n]),$$

(resp. $\mathcal{C}os\mathbf{Top}(\underline{Z}, X_*^{K_*}) \subset \prod_n \mathbf{Top}(\underline{Z}[n], X_*^{K_*}[n])$).

Notice that our results still hold if the category **Top** is replaced by another category with the same formal properties (cocompleteness, existence of a suitable Hochschild-Yoneda functor, and so on). In particular, **Top** can be replaced by the category of simplicial sets $\mathbf{Sets}^{\boldsymbol{\Delta}^{\mathbf{op}}}$.

**Proof.** Notations are as in 1.3. Since $K \otimes_{\boldsymbol{\Delta}} \underline{Z}$ is the colimit in **Top** of the $\underline{Z}[n]$'s over $\mathbf{C} \cong \mathbf{D}$, we have:

$$\begin{aligned}
\mathbf{Top}(K \otimes_{\boldsymbol{\Delta}} \underline{Z}, X) &\cong \\
&\cong \{\phi_\alpha : \underline{Z}[p] \longrightarrow X, \ \alpha \in K_p, \phi_\beta \circ \underline{Z}(f) = \phi_\alpha\} \\
&\cong \{\psi_p : \underline{Z}[p] \longrightarrow X^{K_p} \cong X^K[p] / \ \psi_q \circ \underline{Z}(g) = X^K(g) \circ \psi_p\} \\
&\cong \mathcal{C}os\mathbf{Top}(\underline{Z} X^K).
\end{aligned}$$

where $f \in \mathbf{D}(\alpha, \beta)$, $g \in \boldsymbol{\Delta}([p], [q])\}$ and the $\psi_p$'s are related to the $\phi_\alpha$'s by: $\psi_p(z) := \prod_{\alpha \in K_p} \phi_\alpha(z)$, $z \in \underline{Z}[p]$.

In the pointed case, we can follow the same pattern of proof. We write $*$ for the various base points in $K_*$, $X_*$ and $K_* \otimes_{\boldsymbol{\Delta}} \underline{Z}$ and for the unique element in $\underline{Z}[0]$. We write $\phi_0$ for the map from $\{*\} = \underline{Z}[0]$ to $\{*\} \subset X$.

Then:

**Top**$_*$  $(K_* \otimes_{\mathbf{\Delta}} \underline{Z}, X_*) \cong$

$$\cong \{\phi_\alpha : \underline{Z}[p] \longrightarrow X, \alpha \in K_p|, \begin{cases} \phi_\alpha = \phi_0 \text{ if } \alpha = * \in K_0 \\ \phi_\beta \circ \underline{Z}(f) = \phi_\alpha \end{cases} \}$$
$$\cong \{\psi_p : \underline{Z}[p] \longrightarrow X_*^{K_*}[p]|, \psi_q \circ \underline{Z}(g) = X_*^{K_*}(g) \circ \psi_q\}$$
$$\cong \mathcal{C}\mathbf{osTop}(\underline{Z}, X_*^{K_*}),$$

where $f$ and $g$ are as above.

If we put $\underline{Z} = \underline{\Delta}$ in the formulas in theorem 1 we deduce:

**Corollary 1** *There exist homeomorphisms:*

$$X^{|K|} \cong ||X^K||, \quad \mathbf{Top}_*(|K_*|, X_*) \cong ||X_*^{K_*}||.$$

The first isomorphism is due to Bott-Segal [6, Prop. 5.1]. In the particular case when $K_*$ is the usual reduced pointed simplicial model of $S^1$, the second one is the proposition II,3,5 of Ndombol-Thomas [16].

There is still a homeomorphism $\mathbf{Top}_*(|K_*|, X_*) \cong ||X_*^{K_*}||$ if $K_*$ is not reduced: replace $\underline{Z}$ by $\underline{\Delta}$ and $K_* \otimes_{\mathbf{\Delta}} \underline{Z}$ by $K \otimes_{\mathbf{\Delta}} \underline{\Delta}$, pointed by $(*,*) \in K_0 \times \underline{\Delta}[0]$, in the proof that $\mathcal{C}\mathbf{osTop}(\underline{Z}, X_*^{K_*}) \cong \mathbf{Top}_*(K_* \otimes_{\mathbf{\Delta}} \underline{Z}, X_*)$.

## 2 Cochain algebra of a mapping space.

Let $\Bbbk$ be a field, $K$ a simplicial set and $X$ a space. It should be understood that, from now on and otherwise stated, algebras, tensor products, cochain complexes and so on are defined over $\Bbbk$.

To the cosimplicial space $X^K$ we associate a cochain algebra, written $|NC(X^K)|^*$ and prove that the natural map

$$\phi : |NC(X^K)|^* \to C^*(||X^K||)$$

constructed by Benderski-Gitler [5, Fla.4.5] is a morphism of differential graded algebras (Proposition 2). We then show that, when $K$ is a simplicial finite set weakly equivalent to a finite simplicial set of dimension less or equal to the connectivity of $X$, $\phi$, together with the adjunction isomorphism in the previous section, induces an isomorphism of graded algebras (theorem 2):

$$H^*\left(|NC(X^K)|^*\right) \cong H^*(X^{|K|}).$$

### 2.1 Bicomplex associated to a cosimplicial space.

We write **Vect** for the category of vector spaces over $\Bbbk$ and **Coch** for the category of cochain complexes. Let $S$ be simplicial vector space ($S \in \mathcal{S}\textbf{Vect}$). The *chain complex associated to $S$*, denoted $C_*S$ is the differential graded vector space $(\{S_i\}_{i\geq 0}, \partial)$ with $\partial = d_0 - d_1 + ... + (-1)^n d_n : S_n \to S_{n-1}$. For each degree $n$, $D_n S$ is the subvector space of $S_n$ generated by all the images of the degeneracy operators $s_0, s_1, ..., s_{n-1} : S_{n-1} \to S_n$. It follows from the identities for $d_i \circ s_j$ that $D_*S$ is a subcomplex of $C_*S$. The quotient $N_*S = C_*S/D_*S$ is the normalized chain complex of $S$. The natural projection, $C_*S \to N_*S$ is a chain equivalence, [13, Theorem $VIII - 6.1$]. The cochain complex $C^*S$ (resp. $N^*S$) is the dual (in the graded sense) of $C_*S$ (resp. $(N_*S)$).

Let $X$ be a space and $S = \text{Sing}(X)$ be the simplicial vector space of singular $\Bbbk$-chains over $X$. The cochain complex

$$C^*(X) := C^*\text{Sing}(X) \text{ (resp. } N^*(X) := N^*\text{Sing}(X))$$

is the singular cochain complex of $X$ (resp. the normalized singular cochain complex of $X$). The normalization maps identifies the homology of $C^*(X)$ with the homolology of $N^*(X)$ which are both denoted by $H^*(X)$.

Let $\underline{Z}$ be a cosimplicial topological space. For a fixed $q \geq 0$, $C^q(\underline{Z})$ is a simplicial vector space:

$$(\star^1) \quad C^q(\underline{Z})_p = C^q(\underline{Z}[p]), \ d_i = C^q(d^i), \ s_j = C^q(s^j).$$

Thus $N_*C^q(\underline{Z}) = C_*C^q(\underline{Z})/D_*C^q(\underline{Z})$ is a chain complex and the two induced differentials

$$\begin{aligned}\partial &: N_pC^q(\underline{Z}) \to N_{p-1}C^q(\underline{Z}) \\ \delta &: N_pC^q(\underline{Z}) \to N_pC^{q+1}(\underline{Z})\end{aligned}$$

define a bicomplex, $N_*C^*(\underline{Z})$ whose total complex, denoted $|NC(\underline{Z})|^*$, is a **Z**-graded vector space, with differential $d = \partial - \delta$:

$$|NC(\underline{Z})|^n := \bigoplus_{q-p=n} N_pC^q(\underline{Z}), \quad d : |NC(\underline{Z})|^n \to |NC(\underline{Z})|^{n+1},$$

called the *realization of the simplicial complex* [6, Def.5.2]:

$$[p] \mapsto N_pC^*(\underline{Z}) = \frac{C^*(\underline{Z}[p])}{(C^q(s^j)(C^*(\underline{Z}[p-1])), j = 0, 1...p)}.$$

## 2.2 Cochain algebra of a cosimplicial space.

Recall that, if $X$ is a space, the cup-product [13, 20] provides naturally the complex of singular cochains $C^*(X)$ with the structure of a differential graded algebra.

**Proposition 1** *The complex $|NC(\underline{Z})|^*$ has naturally the structure of a differential graded algebra.*

The product on $|NC(\underline{Z})|^*$ is obtained from the cup product in the $C^*(\underline{Z}[p])$'s and from the shuffles map (see the proof below for explicit formulas). From now on, we view therefore

$$\underline{Z} \mapsto |NC(\underline{Z})|^*$$

as a functor from the category $\mathcal{C}os\mathbf{Top}$ to the category **DA** of differential **Z**-graded algebras. The algebra $|NC(\underline{Z})|^*$ is called *the algebra of cochains on the cosimplicial space $\underline{Z}$*.

**Proof.** Let $A$ and $B$ be two simplicial vector spaces and consider the shuffle map

$$(\star^2) \quad sh : A_p \otimes B_q \to (A \times B)_{p+q}, a \otimes b \mapsto \sum_{\mu,\nu}(-1)^{\epsilon(\mu)}s_\nu a \times s_\mu b, \begin{cases} a & \in A_p \\ b & \in B_q \end{cases}$$

when the sum is taken over the $p+q$ shuffles $\mu_1 < \mu_2 < ... < \mu_p$, $\nu_1 < \nu_2 < ... < \nu_q$, $\mu_i, \nu_j \in \{1, 2, ..., p+q\}$. Here, $\epsilon(\mu) = \sum_{i=1}^p \mu_i - i - 1$ is the signature of the $(p,q)$-shuffle and $s_\mu := s_{\mu_p} \circ s_{\mu_{p-1}} \circ ... \circ s_{\mu_1}$, $s_\nu := s_{\nu_q} \circ s_{\nu_{q-1}} \circ ... \circ s_{\nu_1}$.

Given two integers $r, s \geq 0$ we set $A := C^r(\underline{Z})$, $B := C^s(\underline{Z})$. Then, $A_p = C^r(\underline{Z}[p])$ and $B_q = C^s(\underline{Z}[q])$, and there is a map $m_{p,q}^{r,s}$:

$$C^r(\underline{Z}[p]) \otimes C^s(\underline{Z}[q]) \stackrel{sh}{\to} C^r(\underline{Z}[p+q]) \times C^s(\underline{Z}[p+q]) \stackrel{\cup}{\to} C^{r+s}(\underline{Z}[p+q])$$

given explicitly by the formula:

$$(\star^3) \quad m_{p,q}^{r,s}(x \otimes y) = \sum_{\mu,\nu} (-1)^{\epsilon(\mu)} s_\nu x \cup s_\mu y,$$

with the notations of $(\star^1)$ and $(\star^2)$. Since the shuffle map and the cup product are compatible with the normalization, we get a composite map

$$N_p C^r(\underline{Z}) \otimes N_q C^s(\underline{Z}) \stackrel{sh}{\to} N_{p+q} C^r(\underline{Z}) \times N_{p+q} C^s(\underline{Z}) \stackrel{\cup}{\to} N_{p+q} C^{r+s}(\underline{Z})$$

that induces the product,

$$|NC(\underline{Z})|^{r-p} \otimes |NC(\underline{Z})|^{s-q} \to |NC(\underline{Z})|^{r+s-p-q}, \quad x \otimes y \mapsto x \cdot y.$$

It follows from the naturality property of the shuffle map and the cup product that the product $\cdot$ provides the complex $|NC(\underline{Z})|^*$ with a natural structure of a **Z**-graded differential algebra.

$\diamond$

**Proposition 2** *The natural cochain map [5, Fla.4.5],*

$$\phi : |NC(\underline{Z})|^* \longrightarrow C^*(||\underline{Z}||),$$

*is a homomorphism of differential graded algebras.*

**Proof.** Recall the definition of $\phi$. When $A = Sing(X)$ and $B = Sing(Y)$ the shuffle map, $(\star^2)$ induces a morphism of chain complexes called the Eilenberg-Zilber morphism:

$$EZ : C_p(X) \otimes C_s(Y) \to C_{p+s}(X \times Y).$$

Write $< -, - >$ for the cochain/chain duality and fix $z \in C_p(X)$ so that $\partial z = 0$. The slant product,

$$C^{p+s}(X \times Y) \to C^s(Y), \quad \alpha \mapsto \alpha/z,$$

is defined by the relation

$$< \alpha/z, x > = < \alpha, EZ(z \otimes x) >, \quad x \in C_s(Y).$$

We write $\iota_p \in C_p(\Delta^p)$ for the canonical simplex $id_{\Delta^p}$ and set:

$$\int_{\Delta^p} : C^{p+s}(\Delta^p \times ||\underline{Z}||) \to C^s(||\underline{Z}||), \quad \int_{\Delta^p} \alpha := \alpha/\iota_p.$$

Since $||\underline{Z}|| = \mathcal{C}os\mathbf{Top}(\underline{\Delta}, \underline{Z})$, there are evaluation maps:

$$e_p : \Delta^p \times ||\underline{Z}|| \longrightarrow \underline{Z}[p].$$

Therefore, the composite

$$C_p C^{p+s}(\underline{Z}) = C^{p+s}(\underline{Z}[p]) \xrightarrow{C^*(e_p)} C^{p+s}(\Delta^p \times ||\underline{Z}||) \xrightarrow{\int_{\Delta^p}} C^s(||\underline{Z}||),$$

passes to the normalized bicomplex : $N_p C^{s+p}(\underline{Z}) \to C^s(||\underline{Z}||)$ and induces a cochain complex morphism:

$$\phi : |NC(\underline{Z})|^s \to C^s(||\underline{Z}||).$$

Let $x \in C^{p+r}(\underline{Z}[p])$ and $y \in C^{q+s}(\underline{Z}[q])$. We want to show that:

$$\phi(x \cdot y) = \phi(x) \cup \phi(y),$$

where $\cdot$ stands for the map $m^{*,*}_{p,q} : C^*(\underline{Z}[p]) \otimes C^*(\underline{Z}[q]) \to C^*(\underline{Z}[p+q])$ defined in formula $(\star^3)$.

By definition, $e_p$ is the degree $p$ component of a homomorphism of cosimplicial spaces

$$e : \underline{\Delta} \times ||\underline{Z}|| \longrightarrow \underline{Z},$$

where $\underline{\Delta} \times ||\underline{Z}|| \in \mathcal{C}os\mathbf{Top}$ is cosimplicialy constant in the $||\underline{Z}||$-direction, in the sense that: $\underline{\Delta} \times ||\underline{Z}||[p] := \underline{\Delta}[p] \times ||\underline{Z}||$. Therefore, the cochain algebras morphisms $e_p^* = C^*(e_p)$'s commute with the homomorphisms that are induced at the cochain level by the cosimplicial morphisms acting on $C^*(\underline{Z})$ and $C^*(\underline{\Delta} \times ||\underline{Z}||)$. That is, if $f \in Hom([p],[q])$

$$e_p^* \circ C^*(\underline{Z}(f)) = C^*(f \times id_{||\underline{Z}||}) \circ e_q^*.$$

Since, according to formula $(\star^1)$, the degeneracy operators on $C^{p+s}(\underline{Z})$ (resp. $C^{p+s}(\underline{\Delta} \times ||\underline{Z}||)$) are of the form $s_i^* = C^*(\underline{Z}(s^i))$ (resp. $\sigma_i^* = C^*(s^i \times id_{||z||})$) where $s^i : [p+1] \to [p]$ stands for the canonical $i$-th codegeneracy operator, we get, using $(\star^3)$:

$$\begin{aligned} e^*_{p+q}(x \cdot y) &= \sum_{\mu,\nu}(-1)^{\epsilon(\mu)} \, e^*_{p+q} \circ s^*_\nu(x) \cup e^*_{p+q} \circ s^*_\mu(y) \\ &= \sum_{\mu,\nu}(-1)^{\epsilon(\mu)} \sigma^*_\nu \circ e^*_p(x) \cup \sigma^*_\mu \circ e^*_q(y) \\ &= e^*_p(x) \cdot e^*_q(y). \end{aligned}$$

The last product involved is the product $\cdot$ on $|C^*(\underline{\Delta} \times ||\underline{Z}||)|$, where $\underline{\Delta} \times ||\underline{Z}||$ is viewed as a cosimplicial space.

We are now reduced to prove that:

$$(\star^4) \quad \int_{\Delta^{p_1+p_2}} e_{p_1}^*(x) \cdot e_{p_2}^*(y) = \int_{\Delta^{p_1}} e_{p_1}^*(x) \cup \int_{\Delta^{p_2}} e_{p_2}^*(y).$$

Write $\iota = (\iota_1, \iota_2)$ for the canonical isomorphism:

$$(\Delta^{p_1} \times \Delta^{p_2}) \times ||\underline{Z}|| \times ||\underline{Z}|| \longrightarrow (\Delta^{p_1} \times ||\underline{Z}||) \times (\Delta^{p_2} \times ||\underline{Z}||).$$

We get from the "projection formula" (see [20, Fla.4, p.288]):

$$\int_{\Delta^{p_1}} e_{p_1}^*(x) \cup \int_{\Delta^{p_2}} e_{p_2}^*(y) = \int_{\Delta^{p_1} \times \Delta^{p_2}} \iota_1^*(e_{p_1}^*(x)) \cup \iota_2^*(e_{p_2}^*(y)).$$

The identity $(\star^4)$ follows therefore from the Eilenberg-Zilber formula, which is obtained from the simplicial decomposition of $\Delta^{p_1} \times \Delta^{p_2}$ into copies of the $(p_1+p_2)$-simplex parametrized by the $(p_1, p_2)$-shuffles.

$\diamond$

## 2.3 Convergent cosimplicial mapping spaces.

**Definition 1** *A cosimplicial space $\underline{Z}$ is called convergent if the natural map $\phi : |NC(\underline{Z})|^* \longrightarrow C^*(||\underline{Z}||)$ is a cohomology equivalence.*

If $K$ is a simplicial set and $Y$ is a space, we write $dim(K)$ for the maximal degree of non-degenerate simplices in $K$ and $Conn(Y)$ for the connectivity of $Y$. Recall also that two simplicial sets are weakly equivalent if their geometric realizations are homotopy equivalent and that $K$ is finite if it has only finitely many non degenerate simplices.

**Proposition 3** *Let $K$ be a finite simplicial set and let $Y$ be a space such that $dim(K) \leq Conn(Y)$, then $Y^K$ is convergent.*
*The same statement holds in the pointed case.*

The following is a corollary of propositions 2 and 3:

**Corollary 2** *Under the assumptions of Proposition 3, the natural map of differential graded algebras*

$$\phi : \left|NC(Y^K)\right|^* \mapsto C^*(||Y^K||)$$

*is a quasi-isomorphism.*

Since $|K| \simeq |K'| \Rightarrow Y^{|K|} \simeq Y^{|K'|}$ and since, according to theorem 1 and without restrictions on $K$,

$$||Y^K|| := \mathbf{Costop}(\underline{\triangle}, Y^K) \cong \mathbf{Top}(K \otimes_{\mathbf{\Delta}} \underline{\triangle}, Y) = \mathbf{Top}(|K|, Y) = Y^{|K|}$$

it follows also from Proposition 3 and from the lemma 1 and 2 below that:

**Theorem 2** *Let $K$ be a simplicial finite set (resp. $K_*$ a pointed simplicial finite set) and let $Y$ be a space (resp. $Y_*$ a pointed space). If $K$ (resp. $K_*$) is weakly homotopy equivalent to a finite simplicial set (resp. pointed) of dimension less or equal to $Conn(Y)$, then there is a quasi-isomorphism of differential graded algebras*

$$|NC(Y^K)|^* \cong C^*\left(Y^{|K|}\right) .$$

*The same statement holds in the pointed case.*

In the non-pointed case, Proposition 3 is stated in [5] (Theorem 5.1) under the only hypothesis $dim(K) < Conn(Y)$ and attributed to Bott-Segal. However, Bott and Segal do actually prove a weaker version of the theorem (that we recall below). The proof by induction given in [5] is not correct since the induction hypothesis includes for example the assertion that $Y^K$ is convergent in case $K = \mathcal{S}^n$ is the usual simplicial reduced model of the $n$-sphere (the one with two non-degenerate simplices, one in dimension 0 and one in dimension $n$) and $Y$ is $n+1$-connected. This assertion clearly requires a proof and is far from obvious, even in the rational case (see [18] for a proof of the convergence of $Y^{\mathcal{S}^n}$ in the rational case).

## 2.4 The convergence theorem

Our proof of Proposition 3 is based on the Bott-Segal theorem, which states that proposition 3 is true under the assumption that $K$ is the nerve of a finite open covering of a space or, equivalently, $K$ is a finite *polyhedral simplicial set*. In that case, the theorem follows from the convergence of the Eilenberg-Moore spectral sequence by induction on the number of non-degenerate simplices of $K$ [6]. This proof extends to the pointed case.

Recall here that a simplicial set $K$ is called polyhedral if it is the simplicial set associated to an ordered simplicial complex $L$, i.e. $K = \mathcal{S}L$. In that case each simplex $\sigma \in K_n$ is uniquely determined by the vertices $v_i(\sigma) = d_0 \circ d_1 \circ ... \circ \hat{d_i} \circ ... \circ d_n(\sigma)$. The simplex $\sigma$ is non degenerated iff the vertices are all distinct. Thus each face $d_i\sigma$ is also non degenerated. The geometric realization of a polyhedral set is called a *polyhedron*.

The following two lemmas reduce the general case (pointed or not) to the polyhedral case treated by Bott-Segal. The first one establishes that the geometric realization of any finite dimensional simplicial set is homeomorphic to a polyhedron. More precisely,

**Lemma 1** *Let $K$ be an arbitrary simplicial set (resp. pointed) of dimension $n$. Then, there exist*
*1) $L$, a n-dimensional ordered simplicial complex,*
*2) $\lambda : \mathcal{S}L \to K$ a simplicial set map,*
*3 ) $t : |\mathcal{S}L| \longrightarrow |K|$ a homeomorphism between the geometric realizations, such that $|\lambda|$ is homotopic to $t$. Moreover, the spaces, maps and homotopies are pointed if $K_*$ is pointed.*

The lemma is proven in [12, Theorem 6.1], see also [10, p. 202]. The assertion on $dim(L)$ follows from the construction of $L$, which is a (kind of) double barycentric subdivision. Notice that, if $K$ is finite, then $SL$ is also finite since barycentric subdivision preserves finiteness.

If we write $t^\# : Y^{|K|} \to Y^{|\mathcal{S}L|}$ for the morphism induced by $t$, since $\phi$ is natural, the following diagram is homotopy commutative:

$$\begin{array}{ccc}
|NC(Y^{\mathcal{S}L})|^* & \stackrel{|NC(\lambda^\#)|^*}{\longrightarrow} & |NC(Y^K)|^* \\
\phi_{Y^{\mathcal{S}L}} \downarrow & & \downarrow \phi_{Y^K} \\
C^*(||Y^{\mathcal{S}L}||) & & C^*(||Y^K||) \\
\cong \downarrow & & \downarrow \cong \\
C^*(Y^{|\mathcal{S}L|}) & \stackrel{C^*(t^\#)}{\longrightarrow} & C^*(Y^{|K|})
\end{array}$$

According to the Corollary of Theorem 1 it follows immediately from lemma 1 that

$$||Y^K|| \cong Y^{|K|} \stackrel{t^\#}{\cong} Y^{|\mathcal{S}L|} \cong ||Y^{\mathcal{S}L}||\,.$$

therefore, according to the result of Bott-Segal quoted above, when $K$ is finite, $\phi_{Y^K}$ is a quasi-isomorphism if and only if $|NC(\lambda^\#))|^*$ is a quasi-isomorphism. The same argument shows that, if $\alpha : K \longrightarrow K'$ induces a homotopy equivalence on the realizations and if $\phi_{Y^K}$ (resp. $\phi_{Y^{K'}}$) is a quasi-isomorphism, then $\phi_{Y^{K'}}$ (resp. $\phi_{Y^K}$) is a quasi-isomorphism if and only if $|NC(\alpha^\#)|^*$ is a quasi-isomorphism. In particular, as stated above, Thm. 2 follows therefore from the Lemma 2 below.

Before establishing the lemma, we introduce some notations. Notice the central role played, once again, by left Kan adjunction in the forthcoming constructions and arguments.

Let **Fin** be the full subcategory of **Set** whose objects are the finite sets and recall the existence of the Hochschild-Yoneda functor (1.5) evaluated at a space $Y$:

$$\mathcal{H}(Y) : \mathbf{Set^{op}} \to \mathbf{Top}.$$

The composition of $\mathcal{H}(Y)$ and of the functor of singular cochains with coefficients in $\Bbbk$ (2.2) $C^*(-) : \mathbf{Top^{op}} \to \mathbf{Coch}$ defines, by restriction, a functor:

$$\mathcal{K}(Y) : \mathbf{Fin} \to \mathbf{Coch}.$$

We write $\mathbf{Z[Fin]}$ for the additive category freely generated by **Fin**. Here, additive means that the sets of morphisms are provided with the structure of an abelian group and that the composition of morphisms is bilinear. The objects of $\mathbf{Z[Fin]}$ are the objects of **Fin** and, for $X$ and $Y$ in **Fin**, $Hom_{\mathbf{Z[Fin]}}(X,Y)$ identifies naturally with the free abelian group generated by $Hom_{\mathbf{Fin}}(X,Y)$. Since **Coch** is additive, the functor $\mathcal{K}(Y)$ extends naturally to an additive functor from $\mathbf{Z[Fin]}$ to **Coch**.

Besides, there is a full and faithful 'abelian' Yoneda embedding:

$$\mathcal{Y}_{ab} : \mathbf{Z[Fin]} \longrightarrow \mathbf{Ab}_+^{\mathbf{Z[Fin]^{op}}}, \quad \mathcal{Y}_{ab}[n] := \mathbf{Z[Fin]}(-,[n]).$$

Here and below, the subscript $+$ means that $\mathbf{Ab}_+^{\mathbf{Z[Fin]^{op}}}$ is the category of additive functors from $\mathbf{Z[Fin]^{op}}$ to abelian groups. Since the category **Coch** is cocomplete and since $\mathbf{Z[Fin]}$ is small, we can form the left Kan extension of $\mathcal{K}(Y)$ along $\mathcal{Y}_{ab}$, $Lan_{\mathcal{Y}_{ab}}\mathcal{K}(Y)$ that we prefer once again to write $\mathcal{Y}_{ab} \otimes \mathcal{K}(Y)$:

$$\begin{array}{ccc} \mathbf{Z[Fin]} & \xrightarrow{\mathcal{Y}_{ab}} & \mathbf{Ab}_+^{\mathbf{Z[Fin]^{op}}} \\ {\scriptstyle \mathcal{K}(Y)} \downarrow & \swarrow {\scriptstyle \mathcal{Y}_{ab} \otimes \mathcal{K}(Y)} & \\ \mathbf{Coch} & & \end{array}$$

From the definition of $\mathcal{Y}_{ab}$ we deduce:

$$(\mathcal{Y}_{ab} \otimes \mathcal{K}(Y))(\mathcal{Y}_{ab}[n]) = C^*(\mathcal{H}(Y)[n]) = C^*(Y^{n+1}).$$

We write $\mathcal{Y}_{ab}(K)$ for the simplicial object in $\mathbf{Ab}_+^{\mathbf{Z[Fin]^{op}}}$ defined by:

$$\mathcal{Y}_{ab}(K)_n = \mathcal{Y}_{ab}(K_n) = \mathcal{Y}_{ab}([\#K_n]).$$

Then, $(\mathcal{Y}_{ab} \otimes \mathcal{K}(Y))(\mathcal{Y}_{ab}(K))$ is a cosimplicial cochain complex and there is a cosimplicial cochain complexes isomorphism:

$$C^*(Y^K) \cong (\mathcal{Y}_{ab} \otimes \mathcal{K}(Y))(\mathcal{Y}_{ab}(K)),$$

hence a natural isomorphism between the realizations:

$$|NC(Y^K)|^* \cong |N((\mathcal{Y}_{ab} \otimes \mathcal{K}(Y))(\mathcal{Y}_{ab}(K)))|^*.$$

Recall that if $K \in \mathcal{S}\mathbf{Set}$, then $\Bbbk K$ stands for the associated simplicial vector space and $H_i(K) = H_i(N_*K) = \pi_i(\Bbbk K)$. A map $f \in \mathcal{S}\mathbf{Set}(K, K')$ induces:
    a simplicial map $\Bbbk f \in \mathcal{S}\mathbf{Vect}(\Bbbk K, \Bbbk K')$,
    a cosimplicial map $f^{\#} \in \mathcal{C}os\mathbf{Top}(Y^{K'}, Y^K)$,
    a morphism $|NC(f^{\#})|^* \in \mathbf{DA}(|NC(Y^K)|^*, |NC(Y^{K'})|^*)$.

**Lemma 2** *Let $K$ (resp. $K'$) be a simplicial finite set. If the simplicial map $f : K \to K'$ induces an isomorphism $H_*(K) \to H_*(K')$, then $|NC(f^{\#})|^*$ is a quasi-isomorphism.*

**Proof of Lemma 2.** We refer to [19] and [4] for the generalities on left Kan extensions and abelian categories that we use implicitly in the rest of the proof. The functor $\mathcal{Y}_{ab} \otimes \mathcal{K}(Y)$ between two abelian categories is a left adjoint, and is therefore right exact. Besides, each component $\mathcal{Y}_{ab}(K)_n$ is a projective object in the abelian category $\mathbf{Ab}_+^{\mathbf{Z[Fin]^{op}}}$ (another property of the abelian Yoneda embedding). Therefore, the cohomology of the total complex $|N((\mathcal{Y}_{ab} \otimes \mathcal{K}(Y))(\mathcal{Y}_{ab}(K)))|^*$ depends only of $Y$ and of the isomorphism class of $\mathcal{Y}_{ab}(K)$ in the derived category of the abelian category $\mathbf{Ab}_+^{\mathbf{Z}}[\mathbf{Fin}]$. Since $\mathcal{Y}_{ab}(K)_m[n] = \mathcal{Y}_{ab}(K_m)[n] = \mathbf{Z}[\mathbf{Fin}([n], K_m)]$, $\mathcal{Y}_{ab}(K)_m[n]$ identifies naturally with the group of $m$-dimensional chains on $K^{\times n+1}$. Therefore, $H_*(\mathcal{Y}_{ab}(K))[n] = H_*(K^{\times n+1})$, and $\mathcal{Y}_{ab}(K) \longrightarrow \mathcal{Y}_{ab}(K')$ is a pointwise quasi-isomorphism and a quasi-isomorphism in the category of complexes of $\mathbf{Z}[\mathbf{Fin}]$-modules. The lemma follows in the non-pointed case. The lemma and its proof hold *mutatis mutandis* in the pointed case.

◇

# 3 Finite group actions on mapping spaces.

## 3.1 Cohomological decompositions.

Let $G$ be a finite group and $K$ a simplicial finite set. We assume that there is a (strict, not up to homotopy) action of $G$ on $K$:

$$\mu : G \times K \longrightarrow K.$$

In the pointed case, one should also assume that the base point is a fixed point. Under this hypothesis, the results below apply in the pointed case.

This action induces an action of $G$ on the complex $\mathcal{Y}_{\Bbbk}(K)$ of functors that is associated to $K$. Here, we write $\mathcal{Y}_{\Bbbk}(K)$ for the chain complex associated to the simplicial object in $\mathbf{Vect}_+^{\mathbf{Z}[\mathbf{Fin}]^{op}}$ defined by:

$$\mathcal{Y}_{\Bbbk}(K)_n := \Bbbk \otimes \mathbf{Z}[\mathbf{Fin}](-, K_n).$$

In particular, this construction induces a representation of the group algebra of $G$, also written $\mu$:

$$\mu : \Bbbk[G] \times \mathcal{Y}_{\Bbbk}(K) \longrightarrow \mathcal{Y}_{\Bbbk}(K).$$

The purpose of section 3 is to study what informations on the cohomology of mapping spaces from $|K|$ can be obtained from the properties of this representation. We start by recalling some general properties of group algebras representations and explain how they apply in our setting. We refer for example to [11] for details.

Through this section, we assume that $\Bbbk = \mathbf{F}_p$ and $G$ is of order prime to $p$, or that $\Bbbk$ is a splitting field of characteristic zero for $G$. In these cases, the group algebra $\Bbbk[G]$ is semi-simple and splits as a direct sum of indecomposable 2-sided ideals:

$$\Bbbk[G] \cong I_1 \oplus ... \oplus I_k,$$

where $k$ is the number of conjugacy classes of $G$ and each $I_i$ is a matrix algebra over a division ring. The identity matrix in $I_i$ is a central idempotent in the group algebra that we write $e_i$. The family $(e_i)_{i \in [1,k]}$ is a complete family of orthogonal idempotents in the group algebra, that is:

$$e_i e_j = \delta_i^j e_j, \quad \sum_{i=1}^k e_i = 1.$$

Any module $M$ over the group algebra decomposes canonically as a direct sum:
$$M = \bigoplus_{i=1}^{k} M_i,$$
where $M_i$ is the primary component of $M$ associated to $I_i$ (or, equivalently, associated to the corresponding irreducible representation of $G$, since the set of ideals $I_i$ is canonically in bijection with a full set of non-isomorphic irreducible representations of $G$). In the sequel of the article, we identify $1 \leq i \leq k$ with the corresponding irreducible representation and agree that $i = 1$ corresponds to the trivial representation (the one dimensional representation on which $G$ acts as the identity). The module $M_i$ can be described more explicitly as the image of $M$ under the action of the idempotent $e_i$. These results also apply to actions of $G$ on chain complexes (this follows, for example, from the Schur lemma or from the description of the primary components of representations of $G$ as images under the action of the idempotents) and therefore to the action of $G$ on $\mathcal{Y}_{\Bbbk}(K)$. We get:

**Proposition 4** *The chain complex $\mathcal{Y}_{\Bbbk}(K)$ in $\mathbf{Vect}_+^{\mathbf{Z}[\mathbf{Fin}]^{op}}$ decomposes as a direct sum:*
$$\mathcal{Y}_{\Bbbk}(K) = \bigoplus_{i=1}^{k} \mathcal{Y}_{\Bbbk}^i(K),$$
*where $\mathcal{Y}_{\Bbbk}^i(K) = e_i \mathcal{Y}_{\Bbbk}(K)$ is the primary component of $\mathcal{Y}_{\Bbbk}(K)$ associated to $I_i$.*

**Theorem 3** *Let $Y$ be a $n$-connected space and let $K$ be a simplicial finite set weakly homotopy equivalent to a finite simplicial set of dimension less or equal to $n$. Let $\mu$ be an action of $G$ on $K$ and assume that $\Bbbk$ is as above. Then, $H^*(Y^{|K|})$ splits canonically into primary components under the action of $G$. Moreover, this decomposition can be described as follows:*
$$H^*(Y^{|K|}) \cong \bigoplus_{i=1}^{k} H^*(Y^{|K|})_i \cong \bigoplus_{i=1}^{k} H^*(\mathcal{Y}_{\Bbbk}^i(K) \otimes_{\mathbf{Z}[\mathbf{Fin}]} \mathcal{K}(Y)).$$

Here, $\mathcal{K}(Y)$ is as in subsection 2.4. The action of $G$ on $|K|$ is induced by the action of $G$ on $K$ in the usual way and if $\phi \in Y^{|K|}$ and $g \in G$, $g(\phi) := \phi \circ g^{-1}$. The theorem follows from theorem 2 in the previous section and from proposition 4.

Since, by Yoneda, $\mathcal{Y}_{\Bbbk}(K)$ is a complex of projective and therefore flat objects in $\mathbf{Vect}_+^{\mathbf{Z}[\mathbf{Fin}]^{op}}$, there is a Künneth spectral sequence converging to $H^*(\mathcal{Y}_{\Bbbk}(K) \otimes_{\mathbf{Z}[\mathbf{Fin}]} \mathcal{K}(Y))$ [21]. Theorem 3 implies that this spectral sequence splits into isotypic components. It would be interesting to know if the same argument can be applied to the case of mapping spaces from the spheres and to the action of power maps on them. Let us be more precise. T. Pirashvili has shown that the cohomology groups of mapping spaces from the $n$-spheres decompose into eigenspaces under the action of power maps [18]. His proof relies on natural but technical homological computations in the category of modules over the category of finite pointed sets. In particular, his proof involves a careful study of the above spectral sequence where the category of finite sets is replaced by the category of finite pointed sets and when $K$ is a simplicial model of the $n$-sphere. These results can be explained relatively easily using the algebraic properties of the iterated bar construction [17], but it would be nice to find a direct geometrical proof of Pirashvili's results. Theorem 3 suggest that the pattern of such a proof would read: the action of degree $n$ maps on the spheres translates into an action of the group $\mathbf{Z}$. Pirashvili's decomposition should then be obtained directly from a decomposition into isotypical components of $\mathcal{Y}_{\Bbbk}(K)$ under the action of the group $\mathbf{Z}$. For the purpose of such a proof, one would obviously have to renounce to various of the hypothesis we have worked with ($\mathbf{Z}$ is infinite, the action of $\mathbf{Z}$ should probably be defined only up to homotopy, and so on).

## 3.2 Cohomologically trivial actions on the source.

As a first application of theorem 3, we show in this subsection that, under the hypothesis of theorem 2, a group that acts trivially on the cohomology of the source of a mapping space has also a trivial action on the cohomology of the mapping space. Notice how easily the theorem follows from our previous results. This emphasizes once again the power of Kan adjunction techniques in the study of mapping spaces.

**Theorem 4** *Let $G$, $K$, $\Bbbk$ and $Y$ be as above and assume also that $G$ acts as the identity on the cohomology of $K$. Then, $G$ also acts as the identity on the cohomology of $Y^{|K|}$.*

The theorem follows immediately from the following lemma.

**Lemma 3** *Under the assumptions of theorem 4, we have:*

$$H_*(\mathcal{Y}_{\Bbbk}^i(K)) = 0$$

*whenever $i \neq 1$.*

**Proof.** Recall from the end of section 3 that:

$$H_*(\mathcal{Y}_{\Bbbk}(K)[n]) = H_*(K)^{\otimes n+1}.$$

In particular, $H_*(\mathcal{Y}_{\Bbbk}(K)[n])$ decomposes as:

$$H_*(\mathcal{Y}_{\Bbbk}(K)[n]) = \bigoplus_{i_1,\ldots,i_{n+1}} H_*(K)_{i_1} \otimes \ldots \otimes H_*(K)_{i_{n+1}},$$

where $i_1, \ldots, i_{n+1}$ run from $1$ to $k$, that is, over the set of irreducible $G$-representations. Since, by hypothesis, $H_*(K)_i = 0$ if $i \neq 1$, and since a tensor product of trivial representations is a trivial representation, the lemma follows.

◇

### 3.3 Mapping spaces from Moore spaces.

We conclude this article by computing the splitting of theorem 3 in a simple case. The example explains in particular how the computation can be handled when the isotypic decomposition of $H_*(K)$ is known.

Let $\Bbbk = \mathbf{F}_p$, $p \geq 3$ and $G = \mathbf{Z}/2$. Then, for any pointed simplicial set $K$, $G$ acts on $K \wedge K$ and $|K| \wedge |K|$ by the switching map $\sigma$. In the particular case when $|K|$ is a suspension, this action induces a splitting of $|K| \wedge |K|$ into $\pm 1$ "eigenspaces" [9]:

$$\hat{\sigma} : |K| \wedge |K| \longrightarrow (|K| \wedge |K|)_+ \vee (|K| \wedge |K|)_-.$$

For example, when $|K|$ is the Moore space $M_n = S^n \cup_p e^{n+1}$, this splitting reads:

$$M_n \wedge M_n \cong \Sigma^n M_n \vee \Sigma^{n+1} M_n$$

and the map $M_n \wedge M_n \longrightarrow \Sigma^n M_n$ gives the corresponding $\Sigma$-spectrum $M$ the structure of a ring spectrum [3, 9].

Let now $X$ be a $n+1$-connected pointed space. The mapping space $\mathbf{Top}_*(M_n \wedge M_n, X)$ decomposes up to homotopy as a product:

$$\mathbf{Top}_*(M_n \wedge M_n, X) \cong \mathbf{Top}_*(\Sigma^n M_n, X) \times \mathbf{Top}_*(\Sigma^{n+1} M_n, X)$$

At the cohomological level, this decomposition reads:

$$H^*(\mathbf{Top}_*(M_n \wedge M_n, X)) \cong H^*(\mathbf{Top}_*(\Sigma^n M_n, X)) \otimes H^*(\mathbf{Top}_*(\Sigma^{n+1} M_n, X)).$$

Besides, since $\sigma$ acts trivially on $\Sigma^n M_n$ and as the degree $-1$ map on $\Sigma^{n+1} M_n$, it acts trivially on $H^*(\mathbf{Top}_*(\Sigma^n M_n, X))$ and $H^*(\mathbf{Top}_*(\Sigma^{n+1} M_n, X))$ splits into $\pm 1$ eigenspaces. We compute below the corresponding decomposition on the model for $H^*(\mathbf{Top}_*(M_n \wedge M_n, X))$. We write $K$ for a simplicial model of $M_n$ of dimension $n+1$. Recall that, by theorem 3: $H^*(\mathbf{Top}_*(M_n \wedge M_n, X)) =$

$$H^*(\mathcal{Y}_{\Bbbk}^+(K \wedge K) \otimes_{\mathbf{Z}[\mathbf{Fin}]} \mathcal{K}(X)) \oplus H^*(\mathcal{Y}_{\Bbbk}^-(K \wedge K) \otimes_{\mathbf{Z}[\mathbf{Fin}]} \mathcal{K}(X)).$$

We get two Künneth hyperhomology spectral sequences:

$$\mathbf{Tor}_*^{\mathbf{Z}[\mathbf{Fin}]}(H_*(\mathcal{Y}_{\Bbbk}^+(K \wedge K)), \mathcal{K}(X)) \Rightarrow H^*(\mathcal{Y}_{\Bbbk}^+(K \wedge K) \otimes_{\mathbf{Z}[\mathbf{Fin}]} \mathcal{K}(X)),$$

$$\mathbf{Tor}_*^{\mathbf{Z}[\mathbf{Fin}]}(H_*(\mathcal{Y}_{\Bbbk}^-(K \wedge K)), \mathcal{K}(X)) \Rightarrow H^*(\mathcal{Y}_{\Bbbk}^-(K \wedge K) \otimes_{\mathbf{Z}[\mathbf{Fin}]} \mathcal{K}(X)).$$

By our previous remarks, the first one (resp. the second one) has equivalently $H^*(\mathbf{Top}_*(\Sigma^n M_n, X)) \otimes H^*(\mathbf{Top}_*(\Sigma^{n+1} M_n, X))_+$ as abutment (resp. $H^*(\mathbf{Top}_*(\Sigma^n M_n, X)) \otimes H^*(\mathbf{Top}_*(\Sigma^{n+1} M_n, X))_-$).

We look for example at the $-1$-component. We write $x$ and $y$ for the generators of $H^*(K)$ that appear respectively in the left and right factors in the smash product $K \wedge K$. For a given basis $(x_1, ..., x_4)$ of the reduced cohomology of $K \wedge K$, we write $(x_1^\star, ..., x_4^\star)$ for the dual basis. Then (see e.g. [9]),

$$H_*(K \wedge K)_+ \cong \mathbf{Z}/p \cdot \{(x \otimes y)^\star, (\beta x \otimes y + \beta y \otimes x)^\star\}$$

and

$$H_*(K \wedge K)_- \cong \mathbf{Z}/p \cdot \{(\beta x \otimes y - x \otimes \beta y)^\star, (\beta x \otimes \beta y)^\star\}.$$

Recall that:

$$H_*(\mathcal{Y}_{\Bbbk}(K \wedge K))[n] = \bigoplus_{i_1,...,i_{n+1}} H_*(K \wedge K)_{i_1} \otimes ... \otimes H_*(K \wedge K)_{i_{n+1}},$$

where $i_1, ..., i_{n+1} = \pm 1$. From now on, we write $I$ for the sequence $i_1, ..., i_{n+1}$ and $H_*(K \wedge K)_I$ for $H_*(K \wedge K)_{i_1} \otimes ... \otimes H_*(K \wedge K)_{i_{n+1}}$. We have therefore,

$$H_*(\mathcal{Y}_{\Bbbk}^-(K \wedge K))[n] = \bigoplus_{i_1 i_2 ... i_{n+1} = -1} H_*(K \wedge K)_I.$$

The above basis for $H_*(K \wedge K)_+$ and $H_*(K \wedge K)_-$ provide $H_*(\mathcal{Y}_{\Bbbk}^-(K \wedge K))[n]$ with an explicit basis. The computation of the $E_2$-term of the Künneth hyperhomology spectral sequence follows.

What we would like to point out here, and what the above example shows very explicitly, is that the computation of the isotypic decomposition of the cohomology of a mapping space with a group action on the source relies deeply on the algebra structure of the Grothendieck ring of the group, and not only on the knowledge of the isotypic decomposition of the cohomology of the source. In the above example, this structure is extremely simple since the Grothendieck ring, in that case, is, up to a canonical isomorphism, the group algebra of $\mathbf{Z}/2$. In general, the computation would use the character table of the group.

The following theorem illustrates this idea.

**Theorem 5** *Let $G$, $K$ lk and $Y$ be as in theorem 3. Assume that the cohomology of $K$ splits into a direct sum of isotypical components: $H^*(K) = \bigoplus_{i \in I} H^*(K)_i$, where $I$ is a strict subset of $\{1, 2, ..., k\}$. Then,*

$$H^*(Y^{|K|}) = \bigoplus_{i \in J} H^*(Y^{|K|})_i.$$

*Here, $J$ is a basis for the subalgebra of the algebra of $G$-linear representations (the Grothendieck ring of $G$) generated by $I$.*

The proof follows from our previous computations.